\documentclass[12pt]{amsart}

\usepackage{amsmath, amssymb}
\usepackage{psfig}

\begin{document}

\newtheorem{thrm}{Theorem}
\newtheorem{lmm}[thrm]{Lemma}
\newtheorem{crllry}[thrm]{Corollary}
\newtheorem{cnjctr}[thrm]{Conjecture}

\newcommand{\ignore}[1]{}

\newcommand{\oone}{$\Omega_1$}
\newcommand{\otwo}{$\Omega_2$}
\newcommand{\othr}{$\Omega_3$}
\newcommand{\oen}{$\Omega_n$}

\author{Tobias J. Hagge}

\address{Department of Mathematics, Indiana University, Bloomington, IN 47405}

\email{thagge@indiana.edu}

\date{\today}

\title{Every Reidemeister Move is Needed for Each Knot Type}

\begin{abstract}

We show that every knot type admits a pair of diagrams that cannot be made identical without using Reidemeister \otwo-moves. We also show that our proof is compatible with known results for the other move types, in the sense that every knot type admits a pair of diagrams that cannot be made identical without using all of the move types.

\thanks{The author would like to thank Charles Livingston, Zhenghan Wang, Scott Baldridge, and Noah Salvaterra for their helpful comments, and Vladimir Chernov for pointing out this problem.}
\end{abstract}

\maketitle

\section{Introduction}

Reidemeister proved \cite{Re} that given two diagrams of ambient isotopic links, there is a sequence of transformations on one of the diagrams that gives an explicit isotopy. Each transformation is either a planar isotopy, a cusp move (class \oone), a self-tangency move (class \otwo), or a triple point move (class \othr). We call two diagrams {\em equivalent} if such a sequence exists, and we call the sequence a {\em Reidemeister sequence} for the pair. Sometimes we refer to a Reidemeister sequence without specifying the second diagram in the pair; in this case the second diagram is the result of applying the moves in the sequence. If, for a given diagram pair and $n \in \{1, 2, 3\}$, there is a Reidemeister sequence for the pair that does not contain an \oen-move, we call the pair {\em \oen-independent}. Otherwise the pair is {\em \oen-dependent}.

Since \oone-moves are the only moves that change the winding number of a diagram, it is clear that every link type admits \oone-dependent diagram pairs. Olof-Petter \"{O}stlund \cite{Oe} has shown that every link type admits \othr-dependent diagram pairs, as well as pairs that are simultaneously \oone-dependent and \othr-dependent. In the case of links with at least two components, \otwo-moves are the only moves that change the number of intersections between components. Thus such links admit \otwo-dependent diagram pairs. Vassily Manturov \cite{Ma} has recently shown that a connected sum of any four distinct prime knots admits \otwo-dependent diagram pairs. We consider \otwo-moves in more generality, and construct \otwo-dependent diagram pairs for every knot type. Our approach is similar to Manturov's but was developed independently. We conclude by showing that for each knot type one can construct a diagram pair that is simultaneously \oone-dependent, \otwo-dependent, and \othr-dependent.

\section{Main Theorem}

Briefly, the structure of our main argument is as follows. We give conditions
on a knot diagram that severely limit what one can accomplish without
\otwo-moves. We show that the only transformations possible amount to
replacing the edges in the original diagram with unknotted
$(1,1)$-tangles (by a $(1,1)$-tangle we mean a single stranded tangle). We then show that there is a diagram resulting from a single
\otwo-move that cannot be attained without \otwo-moves. We generalize our result slightly so we can construct \otwo-dependent diagram pairs for every knot type. Finally, we show that our result can be combined with \"{O}stlund's to give a pair of diagrams for each knot type that is simultaneously \oone-dependent, \otwo-dependent, and \othr-dependent.

Let $D$ be a planar knot diagram in general position. A {\em
polygon} $p$ in $D$ is the boundary of a connected component
of the complement $D^c$ of $D$ in the plane. We say $p$ is a {\em
0-gon} if $D$ contains no crossing points. We say $p$ is an {\em
n-gon} if, when all crossing points of $D$ that lie on $p$ are
removed, the remainder consists of $n$ connected components
homeomorphic to an open interval. The points so removed are called the
{\em vertices} of $p$, and the connected components are the {\em
edges}. Note that $p$ can have fewer vertices than edges, and for our purposes a 0-gon has zero edges.

\begin{thrm}
\label{mainthrm}
Let $D$ be a diagram for which the following hold:

\begin{enumerate}

\item There are no 0-gons, 1-gons, or 2-gons,

\item The first \oen-move cannot be an \othr-move.

\end{enumerate}

Then any Reidemeister sequence not containing an \otwo-move does nothing more up to planar isotopy than replace the edges of $D$ with unknotted $(1,1)$-tangles. Furthermore, one can apply a single \otwo-move that crosses two distinct edges of $D$, giving a diagram $D'$ such that the pair $\{D, D'\}$ is \otwo-dependent.
 
\end{thrm}

\begin{figure}
\centerline{\psfig{figure=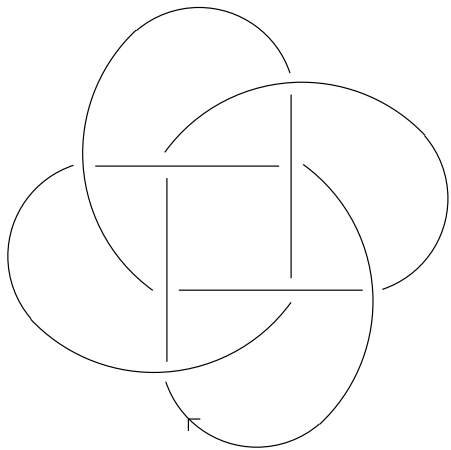}}
\caption{A simple diagram satisfying the preconditions of Theorem \ref{mainthrm}.}\label{turks}
\end{figure}

\begin{figure}
\centerline{\psfig{figure=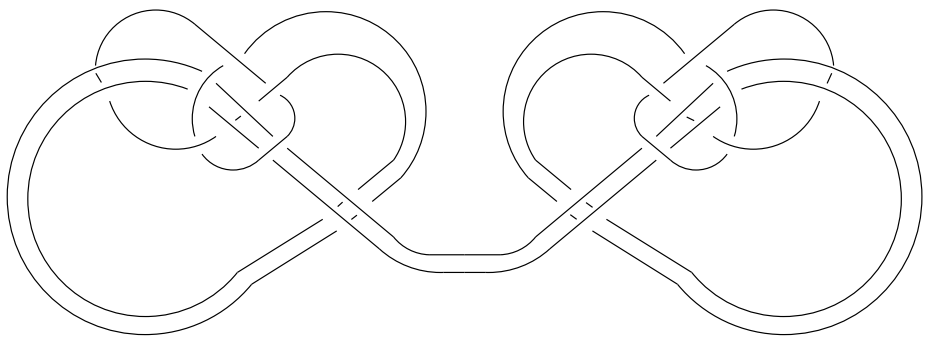}}
\caption{A diagram of the unknot satisfying the preconditions of Theorem \ref{mainthrm}.}\label{unknot}
\end{figure}

Figure~\ref{turks} satisfies the above preconditions. Figure~\ref{unknot} gives
an unknotted example. There are other examples, such as alternating
diagrams with no 0-gons, 1-gons, or 2-gons.

\begin{proof}
Suppose from now on that we do not make any \otwo-moves. Let $\{k_i\}$ be the set of crossing points in $D$. In order for one of the $k_i$ to take part in an \oone-move (in the sense that $k_i$ appears in at least one of the pictures that locally describe the move), it must be the vertex of a 1-gon. In order for a $k_i$ to take part in an \othr-move, it must be one of the vertices of a 3-gon. Since $D$ contains no 1-gons, and none of its 3-gons admit an \othr-move, no $k_i$ can take part in the first \oen-move.

Suppose then that we have made a sequence of moves, and no $k_i$ has taken part in an \oone-move or an \othr-move. Then our diagram is isotopic to the original diagram with each of the edges replaced by an unknotted $(1,1)$-tangle (Figure~\ref{local} provides an illustration, where each dotted box contains a $(1,1)$-tangle). The $(1,1)$-tangles cannot intersect because any \oone-move in which no $k_i$ takes part is just a kink on a single $(1,1)$-tangle, and any \othr-move in which no $k_i$ takes part cannot cause the intersection of two $(1,1)$-tangles that did not intersect before the move.

\begin{figure}
\centerline{\psfig{figure=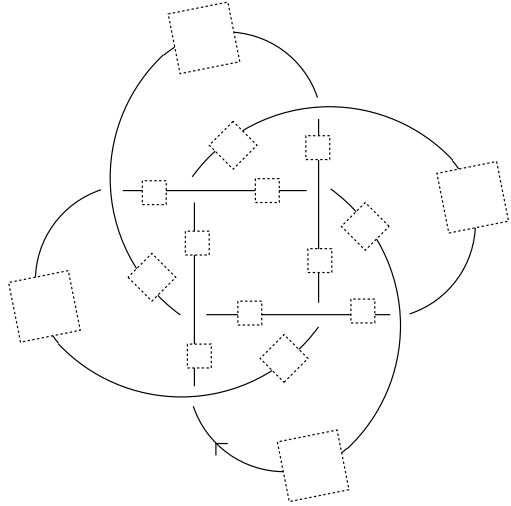}}
\caption{Figure~\ref{turks}, with edges replaced by unknotted $(1,1)$-tangles. Each dotted box represents an unknotted $(1,1)$-tangle.}\label{local}
\end{figure}

Now, it is still impossible for one of the $k_i$ to take part in an \oone-move or an \othr-move. For, given a $k_i$, every polygon $p$ that contains $k_i$ as a vertex looks like one of the polygons in $D$ that contain $k_i$, possibly with some extra edges due to the $(1,1)$-tangles. Thus, $p$ has at least three edges, and $p$ has exactly three edges only when it is one of the 3-gons in $D$ (here we need that there are no 2-gons in $D$). Since none of the 3-gons in $D$ admit \othr-moves, no $k_i$ can take part in an \oone-move or an \othr-move.

Thus, up to isotopy, if \otwo-moves are not allowed then a sequence of moves on $D$ will fix the $k_i$ and replace the edges of $D$ with unknotted $(1,1)$-tangles. This gives the first part of the theorem.

For the second part of the theorem, we create a new diagram $D'$ by performing a single \otwo-move that crosses two distinct edges of $D$. Figure~\ref{turkstwo} provides an example. Let $E$ be a diagram such that $(D,E)$ is an \otwo-independent diagram pair. There are easy ways to show that $D'$ and $E$ cannot be isotopic. For instance, one could show that if $D'$ and $E$ have the same number of crossings, $E$ must contain a 1-gon, while $D'$ contains no 1-gons. Instead we present an argument using Gauss diagrams that straightforwardly generalizes to our full result.

\begin{figure}
\centerline{\psfig{figure=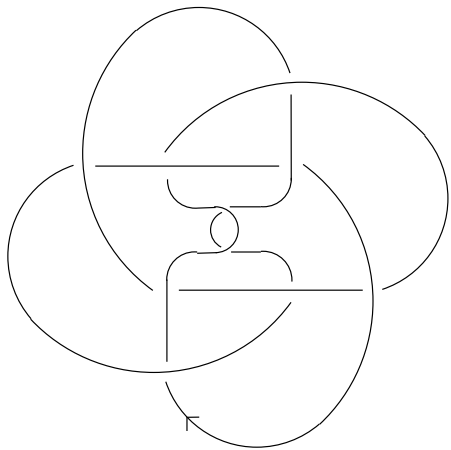}}
\caption{This diagram cannot be derived from Figure~\ref{turks} without an \otwo-move.}\label{turkstwo}
\end{figure}

Every oriented knot diagram $B$ is given by a smooth immersion $\phi$ from the oriented circle $S^1$ to the plane $\mathbb R^2$ with decorated crossing points. This immersion is unique up to orientation preserving self-diffeomorphisms of $S^1$. The map $\phi$ is one to one except at crossing points, where it is two to one. The {\em Gauss diagram} $G_B$ for $B$ is constructed from $S^1$ by drawing a signed arrow between the two elements of $\phi^{-1}(k)$ for each crossing point $k$ of $B$. Each arrow points toward the over strand of the crossing. The sign of each arrow is the sign of that crossing, either $+1$ or $-1$ according to the standard convention. Figures \ref{turksgauss} and \ref{eight} give examples.

\begin{figure}
\centerline{\psfig{figure=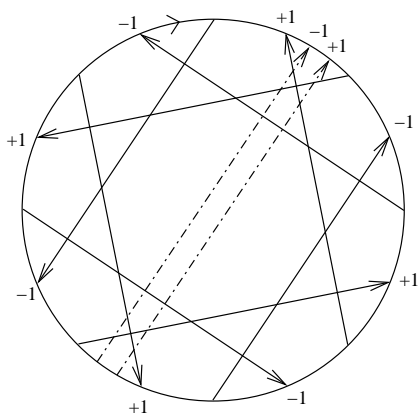}}
\caption{The Gauss diagram for Figure~\ref{turkstwo}. Removing the dashed arrows gives the Gauss diagram for Figure~\ref{turks}.\label{turksgauss}}
\end{figure}

Let $D$ be given an orientation. This gives an orientation for $E$. Consider the Gauss diagrams $G_D$ and $G_E$ of $D$ and $E$. One can see that $G_E$ is just $G_D$ with some extra arrows representing the $(1,1)$-tangles. It is easy to verify that none of the extra arrows intersect the arrows of $G_D$.

The \otwo-move that takes $D$ to $D'$ adds two arrows $a_1$ and $a_2$ to $G_D$, giving the Gauss diagram $G_{D'}$. The heads of $a_1$ and $a_2$ are adjacent on $S^1$, as are the tails. Thus, any arrow that intersects $a_1$ or $a_2$ must intersect both $a_1$ and $a_2$. One can easily show that both $a_1$ and $a_2$ must intersect at least one of the arrows in $G_D$.

Any Gauss diagram containing a pair of arrows with adjacent heads and tails represents a knot diagram that contains a 2-gon. Since $D$ contains no 2-gons, any copy of $G_D$ in $G_{D'}$ intersecting no other arrows cannot contain $a_1$ or $a_2$. Thus the addition of the arrows $a_1$ and $a_2$ causes $G_D$ to intersect other arrows without creating any new copies of $G_D$ that don't intersect other arrows. This reduces the number of copies of $G_D$ intersecting no other arrows from one to zero. Since $E$ has at least one copy of $G_D$ intersecting no other arrows, $D'$ and $E$ are not isotopic. This proves the theorem.
\end{proof}

\begin{thrm}
Every knot type admits an \otwo-dependent diagram pair.\label{bettermain}
\end{thrm}

\begin{proof}
Figure~\ref{unknot} is a diagram of the unknot and so Theorem~\ref{mainthrm} gives us an \otwo-dependent unknotted diagram pair. If we remove a small closed line segment from one of the edges of a diagram $A$, we are left with a $(1,1)$-tangle $T$ having the same knot type as $A$. If $B$ is another knot diagram, then $A\#B$ is obtained by replacing an edge in $B$ with a $(1,1)$-tangle planar isotopic to $T$.

Let $D$ be the diagram in Figure~\ref{unknot}. Suppose we replace one or more of the edges of $D$ with arbitrary $(1,1)$-tangles and get a diagram $F$. Then by the same argument as in Theorem~\ref{mainthrm}, the only thing we can do to $F$ without using \otwo-moves is to replace these $(1,1)$-tangles with other $(1,1)$-tangles of the same knot type. Now, $G_F$ must contain at least one copy of $G_D$ that intersects no other arrows. It may contain more (for instance if $F=D\#D$). However, none of these copies can be altered or intersected with other arrows by a Reidemeister sequence not containing \otwo-moves. We can make a single \otwo-move that crosses two distinct edges of $F$ to get a diagram $F'$ such that the created arrows on $G_{F'}$ intersect a copy of $G_D$ that intersected nothing else in $G_F$. This reduces the number of copies of $G_D$ that intersect no other arrows, just as in Theorem~\ref{mainthrm}. Thus $(F, F')$ is an \otwo-dependent diagram pair.
\end{proof}

We wish to construct diagram pairs that are simultaneously \oone-dependent, \otwo-dependent, and \othr-dependent. In order to do this we briefly summarize a portion of the proof given by \"{O}stlund in \cite{Oe}, enough to prove the existence of \othr-dependent diagram pairs. The reader should see \cite{Oe} for details.

\begin{thrm}
Every knot type admits an \oone-dependent, \otwo-dependent, \othr-dependent diagram pair.
\end{thrm}
\begin{proof}

\"{O}stlund's proof counts the signed number of instances of the Gauss subdiagram given in Figure~\ref{hash}. The sign of each subdiagram is given by the product of the signs of the crossings in the subdiagram. \"{O}stlund shows that this count is invariant under \oone-moves and \otwo-moves, but can vary under \othr-moves. \"{O}stlund uses this count to prove that for every knot type there is a diagram pair that is simultaneously \oone-dependent and \othr-dependent. As an example, the Gauss diagram in Figure~\ref{eight} represents a figure eight knot. It contains one copy of Figure~\ref{hash} as a subdiagram with sign $1^2(-1)^2=1$. Its mirror image has the same Gauss diagram with the arrows and signs (but not the orientation of $S^1$) reversed, and this does not contain any copies of Figure~\ref{hash} as a subdiagram. Figure~\ref{eight} and its mirror image have different winding numbers, so the pair is \oone-dependent and \othr-dependent.

\begin{figure}
\centerline{\psfig{figure=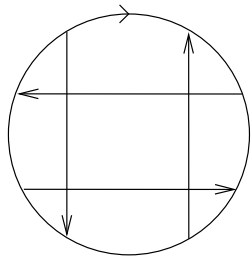}}
\caption{The signed sum of subdiagrams of this form is invariant under \oone-moves and \otwo-moves, but not \othr-moves.}\label{hash}
\end{figure}

\begin{figure}
\centerline{\psfig{figure=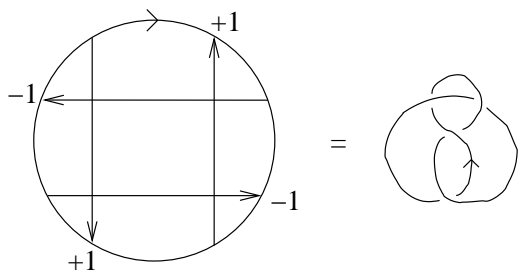}}
\caption{A Gauss diagram for a figure eight knot.}\label{eight}
\end{figure}

\"{O}stlund's count is also additive under connected sum of diagrams. Thus, if $(D, D')$ is an \otwo-dependent diagram pair constructed according to Theorem~\ref{bettermain}, where $D$ and $D'$ differ only by an \otwo-move, and $(E, E')$ is one of \"{O}stlund's \oone-dependent, \othr-dependent diagram pairs, then the diagrams $D \# E$ and $D' \# E'$ will form an \oone-dependent, \otwo-dependent, \othr-dependent diagram pair.
\end{proof}

\begin{figure}
\centerline{\psfig{figure=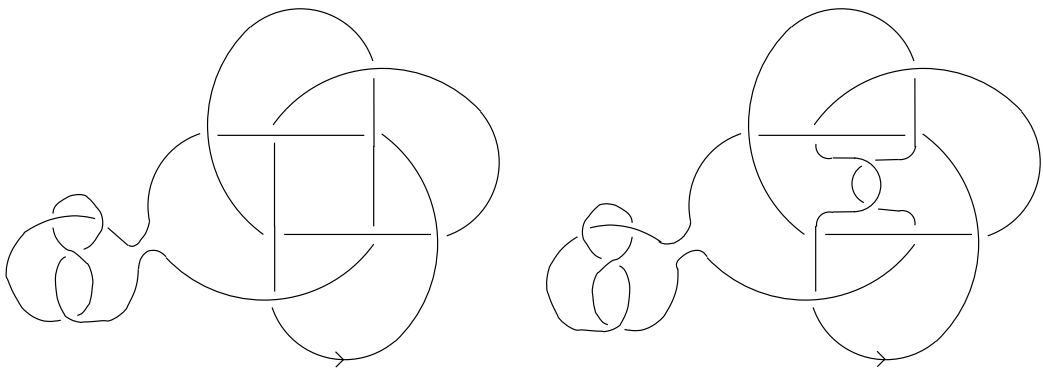}}
\caption{These diagrams have different winding numbers. The pair is \oone-dependent, \otwo-dependent, and \othr-dependent.}\label{dependent}
\end{figure}

 Figure~\ref{dependent} gives an example.

\bibliographystyle{plain}
\bibliography{references}
\end{document}